\documentclass{amsart}

\usepackage{graphicx}

\usepackage{}

\newtheorem{theorem}{Theorem}[section]

\newtheorem{corollary}[theorem]{Corollary}
\newtheorem{conjecture}[theorem]{Conjecture}

\theoremstyle{definition}

\theoremstyle{remark}

\numberwithin{equation}{section}

\newcommand{\cF}{{\mathcal{F}}}

        \newcommand{\field}[1]{{\mathbb{#1}}}
        \newcommand{\NN}{\field{N}}
        \newcommand{\ZZ}{\field{Z}}
        
        \newcommand{\RR}{\field{R}}

\newcommand{\Dom}{\mbox{\rm Dom}}

\newcommand{\Tr}{\mbox{\rm Tr}}

\begin{document}

\title[Schr\"odinger operators with magnetic wells]{Semiclassical analysis of Schr\"odinger operators with magnetic wells}


\author{Bernard Helffer}
\address{D\'epartement de Math\'ematiques, B\^atiment 425, Univ
Paris-Sud et CNRS, F-91405 Orsay C\'edex, France}
\email{Bernard.Helffer@math.u-psud.fr}

\author{Yuri A. Kordyukov}
\address{Institute of Mathematics, Russian Academy of Sciences, 112 Chernyshevsky
str. 450077 Ufa, Russia} \email{yurikor@matem.anrb.ru}
\thanks{Y.K. was partially supported by the Russian Foundation of Basic
Research (grant 06-01-00208).}


\date{}

\begin{abstract}
We give a survey of some results, mainly  obtained  by the authors and
their collaborators, on spectral properties of the magnetic
Schr\"odinger operators in the semiclassical limit. We focus our
discussion on asymptotic behavior of the individual eigenvalues for
operators on closed manifolds and existence of gaps in intervals
close to the bottom of the spectrum of periodic operators.
\end{abstract}

\maketitle

\section{Preliminaries}
\subsection{The magnetic Schr\"odinger operators}
Let $(M,g)$ be an oriented Riemannian manifold of dimension $n\geq
2$. Let $\bf B$ be a real-valued closed $C^\infty$ 2-form on $M$.
Assume that $\bf B$ is exact and choose a real-valued $C^\infty$
1-form $\bf A$ on $M$ such that $d{\bf A} = \bf B$.

Thus, one has a natural mapping
\[
u\mapsto i\,du+{\bf A}u
\]
from $C^\infty_c(M)$ to the space $\Omega^1_c(M)$ of smooth,
compactly supported one-forms on $M$. The Riemannian metric allows
to define scalar products in these spaces and consider the adjoint
operator
\[
(i\,d+{\bf A})^* : \Omega^1_c(M)\to C^\infty_c(M).
\]
A Schr\"odinger operator with magnetic potential $\bf A$ is defined
by the formula
\[
H_{\bf A} = (i\,d+{\bf A})^* (i\,d+{\bf A}).
\]

From the geometric point of view, we may regard $\bf A$ as a
connection one form of a Hermitian connection on the trivial line
bundle $\mathcal L$ over $M$, defining the covariant derivative
$\nabla_{\bf A} = d-i{\bf A}$. The curvature of this connection is
$-i\bf B$. Then the operator $H_{\bf A}$ coincides with the
covariant (or Bochner) Laplacian:
\[
H_{\bf A}=\nabla_{\bf A}^*\nabla_{\bf A}.
\]

Choose local coordinates $X=(X_1,\ldots,X_n)$ on $M$. Write the
1-form $\bf A$ in the local coordinates as
\[
{\bf A}= \sum_{j=1}^nA_j(X)\,dX_j,
\]
the matrix of the Riemannian metric $g$ as
\[
g(X)=(g_{j\ell}(X))_{1\leq j,\ell\leq n}
\]
and its inverse as
\[
g(X)^{-1}=(g^{j\ell}(X))_{1\leq j,\ell\leq n}.
\]
Denote $|g(X)|=\det(g(X))$. Then the magnetic field $\bf B$ is given
by the following formula
\[
{\bf B}=\sum_{j<k}B_{jk}\,dX_j\wedge dX_k, \quad
B_{jk}=\frac{\partial A_k}{\partial X_j}-\frac{\partial
A_j}{\partial X_k}.
\]
Moreover, the operator $H_{\bf A}$ has the form
\[
H_{\bf A}=\frac{1}{\sqrt{|g(X)|}}\sum_{1\leq j,\ell\leq n}\left(i
\frac{\partial}{\partial X_j}+A_j(X)\right) \left[\sqrt{|g(X)|}
g^{j\ell}(X) \left(i \frac{\partial}{\partial
X_\ell}+A_\ell(X)\right)\right].
\]

When $n=2$, the magnetic two-form $\bf B$ is a volume form on $M$
and therefore can be identified with the function $b\in C^\infty(M)$
given by
\[
{\bf B}=b\,dx_g,
\]
where $dx_g$ denotes the Riemannian volume form $M$ associated with
$g$.

When $n=3$, the magnetic two-form $\bf B$ can be identified with a
magnetic vector field $\vec{b}$ by the Hodge star-operator. If $M$
is the Euclidean space $\RR^3$, we have
\[
\vec{b}=(b_1,b_2,b_3) = \operatorname{curl} {\bf
A}=(B_{23},-B_{13},B_{12})
\]
with the usual definition of $\operatorname{curl}$.

We will consider the magnetic Schr\"odinger operator $H_{\bf A}$ as
an unbounded operator in the Hilbert space $L^2(M)$. We will discuss
two cases:
\begin{itemize}
  \item $M$ is a compact manifold, possibly with boundary;
  \item $M$ is a noncompact oriented manifold equipped with a properly
discontinuous action of a finitely generated, discrete group
$\Gamma$ such that $M/\Gamma$ is compact.
\end{itemize}

In the first case, if $M$ has non-empty boundary, we will assume that the operator
$H^h$ satisfies the Dirichlet boundary conditions. Moreover, we will
only consider the case when the potential wells defined by the
magnetic field lie in the interior of $M$.
A closely related case is the case $M=\mathbb R^2$ under the
assumption that the potential wells defined by the magnetic field
lie in a compact subset of $\mathbb R^2$ and that
$\liminf_{|x|\rightarrow
+\infty}
 b(x) > \inf b$.

In the second case, we will assume that
    $M$ is complete and
$H^1(M, \RR) = 0$, i.e. any closed $1$-form on $M$ is exact.
Moreover, the metric $g$ and the magnetic 2-form $\bf B$ are
supposed to be $\Gamma$-invariant (but $\bf A$, in general, is not
$\Gamma$-invariant). Moreover, we will assume that the magnetic
field has a periodic set of compact potential wells (see
Section~\ref{s:periodic} for a precise definition).

In both cases, if $M$ is without boundary (this is always
true in the second case), the operator $H_{\bf A}$ is essentially
self-adjoint with domain $C^\infty_c(M)$. In the case when $M$ has
non-empty boundary, we will consider the self-adjoint operator
obtained as the Friedrichs extension of the operator $H_{\bf A}$
with domain $C^\infty_c(M)$ (the Dirichlet realization).
We refer the reader to the book
\cite{Fournais-Helffer:book} (and the references therein),
for the description of the  spectral properties of  the Neumann realization of
 a magnetic Schr\"odinger operator
on a compact manifold with boundary and their applications to problems in
superconductivity and liquid crystals. We also refer the reader to
the surveys \cite{CFKS,Helffer-Kuroda,Helffer-ECM,Mohamed-Raikov}
for the presentation of general results concerning the Schr\"odinger
operator with magnetic fields.

We will discuss spectral properties of the magnetic Schr\"odinger
operator in the semiclassical limit. So we consider the operator
$H^h$, depending on a semiclassical parameter $h>0$, defined as
\[
H^h = (ih\,d+{\bf A})^* (ih\,d+{\bf A}).
\]
The operators $H^h$ and $H_{\bf A}$ are related by the formula
\[
H^h=h^2(\,d-ih^{-1} {\bf A})^* (\,d-ih^{-1}{\bf A})= h^2 H_{h^{-1}
\bf A}.
\]
This formula shows, in particular, that the semiclassical limit
$h\to 0$ is clearly equivalent to the large magnetic field limit.

\subsection{Magnetic wells}
For any $x\in M$, denote by $B(x)$ the linear operator on the
tangent space $T_x{ M}$ associated with the 2-form $\bf B$:
\[
g_x(B(x)u,v)={\bf B}_x(u,v),\quad u,v\in T_x{ M}.
\]

In local coordinates $X=(X_1,\ldots,X_n)$, the matrix
$(b^\alpha_\beta(X))_{\alpha,\beta=1,\ldots,n}$ of $B(X)$ is given
by
\[
b^\alpha_\beta(X)=\sum_{j=1}^n B_{\beta j}(X)g^{j\alpha}(X).
\]

It is easy to check that $B$ is skew-adjoint with respect to $g$,
and therefore for each $x\in M$ the non-zero eigenvalues of $B(x)$
can be written as $\pm i\lambda_j(x)$, where $\lambda_j(x)>0$,
$j=1,2,\ldots, d$. Introduce the function (the intensity of the
magnetic field)
\[
{\Tr}^+ (B(x))=\sum_{j=1}^d \lambda_j(x)=\frac{1}{2}\Tr([B^*(x)\cdot
B(x)]^{1/2}).
\]

We will also use the trace norm of $B(x)$:
\[
|B(x)|=[\Tr (B^*(x)\cdot B(x))]^{1/2}\,.
\]
It coincides with the norm of $B(x)$ with respect to the Riemannian
metric on the space of tensors of type $(1,1)$ on $T_xM$ induced by
the Riemannian metric $g$ on $M$. In local coordinates
$X=(X_1,\ldots,X_n)$, we have
\[
|B(X)|=\left(\sum_{i,j,k,\ell} g^{ij}(X)g_{k\ell}(X) b_i^k(X)
b_j^\ell(X)\right)^{1/2}\,.
\]
When $n=2$, then
\[
{\Tr}^+ (B(x))=|b(x)|\ \text{and}\ |B(x)|=\sqrt{2}\,|b(x)|\,.
\]
When $n=3$, then
\[
{\Tr}^+ (B(x))=|\vec{b}(x)|\ \text{and}\
|B(x)|=\sqrt{2}\,|\vec{b}(x)|\,.
\]

Remark that the function $|B(x)|^2$ is clearly $C^\infty$, whereas
the function ${\Tr}^+ B$ is only continuous (more precisely, it is
locally H\"older of order $1/2n$ (see \cite{Helffer-Robert84} and
references therein)). It turns out that in many spectral problems
the function $x\mapsto h\cdot {\Tr}^+ (B(x))$ can be considered as a
magnetic potential, that is, as a magnetic analog of the electric
potential $V$ in a Schr\"odinger operator $-h^2\Delta+V$. This leads
us to introduce the notion of magnetic well as follows.

Let $b_0$ be the minimal intensity of the magnetic field
\[
b_0=\min \{{\Tr}^+ (B(x))\, :\, x\in { M}\}.
\]
Consider the zero set of ${\Tr}^+ (B(x))-b_0$
\[
U=\{x\in M : {\Tr}^+ (B(x))=b_0\}\,.
\]
A magnetic well (attached to the given energy $hb_0$) is by
definition a connected component of $U$. If $M$ is compact and has non-empty boundary, we will always
assume that $U$ is included in the interior of $M$.

\subsection{Rough estimates for the lowest eigenvalue}
Assume that $M$ is a compact manifold. Denote by $\lambda_0(H^h)$
the bottom of the spectrum of the operator $H^h$ in $L^2(M)$.

\begin{theorem}[\cite{HM},
Theorem 2.2]
For any $\mu\in \operatorname{Im} {\Tr}^+ B$, there exists $C>0$
 and $h_0>0$
such that, for any $h\in (0,h_0]$\,,
\[
(-C\,h^{4/3}+h\mu,h\mu+C\,h^{4/3})\neq \emptyset\,.
\]
Moreover, there exists $C>0$ such that
\[
-C \,h^{5/4} \leq \lambda_0(H^h)-hb_0 \leq C\, h^{4/3}.
\]
\end{theorem}

The last result can be improved if the rank of ${\bf B}$ is
constant. This can be seen as a form of the Melin-H\"ormander
inequality. Using the techniques developed in
\cite{Helffer-Robert84} one can indeed get the existence of $C>0$
 and $h_0>0$ such that,  for any $h\in (0,h_0]$\,,
\[
-C \,h^2 \leq \lambda_0(H^h)-hb_0\,.
\]

Remark that if $n=2$ and $M$ is without boundary
 then we necessarily have $b_0=0$, since
\[
\int_M b(x)dx_g=\int_Md{\bf A}=0\,.
\]
If we suppose that $M$ has non-empty boundary, the operator $H^h$
satisfies the Dirichlet boundary conditions and $b_0>0$, it was
observed by many authors
\cite{Mont,Mall,Ue} (as the immediate consequence of the
Weitzenb\"ock-Bochner type identity and the positivity of the square
of a suitable Dirac operator) that
\[
\inf \sigma(H_{\bf A})\geq b_0\,,
\]
where $\sigma(H_{\bf A})$ denotes the spectrum of the operator
$H_{\bf A}$ in $L^2(M)$ and, as a consequence, that, for any
$h>0\,$,
\[
\lambda_0(H^h)\geq hb_0\,.
\]
In the case $M=\mathbb R^2$,
 this estimate follows from the formula
\[
b(x) = -i [D_{x_1} -A_1, D_{x_2}-A_2]\,,
\]
where, as usual, $D_{x_k} =\frac1i\frac{\partial }{\partial x_k}$,
$k=1,2$, which implies (after an integration by parts) that
\[
\int b(x) |u(x)|^2\,dx \leq  \|(D_{x_1} - A_1) u\|^2 + \|(D_{x_2} -
A_2) u\|^2\;.
\]

If $b_0=0$, one prove a more precise estimate for $\lambda_0(H^h)$
in the case when the magnetic wells are regular submanifolds. Denote
by $d(x,y)$ the geodesic distance between $x$ and $y$.

\begin{theorem}[\cite{HM}, Theorem 2.4]
Let us assume that $b_0=0$ and that $U$ is a $C^\infty$ compact submanifold
of $M$  included in the interior in $M$.
If there exist $k\in\ZZ_+$, $C_1$ and $C_2>0$ such that if
$d(x,U)<C_2$
\[
C_1^{-1}\,d(x,U)^k\leq |B(x)|  \leq C_1\, d(x,U)^k\,,
\]
then one can find $h_0$ and $C>0$ such that, for any $h\in (0,h_0]\,$,
\[
C^{-1}\, h^{(2k+2)/(k+2)} \leq \lambda_0(H^h) \leq C\, h^{(2k+2)/(k+2)}\,.
\]
\end{theorem}
\

\section{Discrete wells}
In this section, we continue to assume that $M$ is compact. Denoting
by $\lambda_0(H^h)\leq \lambda_1(H^h)\leq \lambda_2(H^h)\leq \ldots$
the eigenvalues of the operator $H^h$ in $L^2(M)$, we will consider
the case when the magnetic wells are points.

\subsection{The case $b_0=0$}
Let us assume that $b_0=0$, and, for some integer $k>0$, if
$B(x_0)=0$, then
    $x_0$ belongs to the interior of $M$ and
there exists a positive constant $C$ such that for all $x$ in some
neighborhood of $x_0$ the estimate holds:
\[
C^{-1}\, d(x,x_0)^k\leq {\Tr}^+ (B(x))  \leq C\, d(x, x_0)^k.
\]

In this case, the important role is played by a differential
operator $K^h_{\bar{x}_0}$ in $\RR^n\,$, which is in some sense an
approximation to the operator $H^h$ near $x_0\,$. Recall its
definition \cite{HM}.

Let $\bar x_0$ be a zero of $B$. Choose local coordinates $f: U(\bar
x_0)\to \RR^n$ on $M$, defined in a sufficiently small neighborhood
$U(\bar x_0)$ of $\bar x_0$. Suppose that $f(\bar x_0)=0$, and the
image $f(U(\bar x_0))$ is a ball $B(0,r)$ in $\RR^n$ centered at the
origin.

Write the $2$-form $\bf B$ in the local coordinates as
\[
{\bf B}(X)=\sum_{1\leq l<m\leq n} b_{lm}(X)\,dX_l\wedge dX_m\,, \quad
X=(X_1,\ldots,X_n)\in B(0,r).
\]
Let ${\bf B}^0$ be the closed 2-form in $\RR^n$ with polynomial
components defined by the formula
\[
{\bf B}^0(X)=\sum_{1\leq l<m\leq
n}\sum_{|\alpha|=k}\frac{X^\alpha}{\alpha !}\frac{\partial^\alpha
b_{lm}}{\partial X^\alpha}(0)\,dX_l\wedge dX_m\,, \quad X\in\RR^n.
\]
One can find a 1-form ${\bf A}^{0}$ on $\RR^n$ with polynomial
components such that
\[
d{\bf A}^0(X) ={\bf B}^0(X)\,, \quad X\in\RR^n.
\]

Let $K^h_{\bar{x}_0}$ be a self-adjoint differential operator in
$L^2(\RR^n)$ with polynomial coefficients given by the formula
\[
K_{\bar{x}_0}^h =  (i h\,d+{\bf A}^0)^* (i h\,d+{\bf A}^0)\,,
\]
where the adjoints are taken with respect to the Hilbert structure
in $L^2(\RR^n)$ given by the flat Riemannian metric $(g_{lm}(0))$ in
$\RR^n$. If ${\bf A}^0$ is written as
\[
{\bf A}^0=A^0_{1}\, dX_1+\ldots+ A^0_{n}\,dX_n\,,
\]
then $K^h_{\bar{x}_0}$ is given by the formula
\[
K_{\bar{x}_0}^h=\sum_{1\leq l,m\leq n} g^{lm}(0) \left(i h
\frac{\partial}{\partial X_l}+A^0_{l}(X)\right)\left(i h
\frac{\partial}{\partial X_m}+A^0_{m}(X)\right).
\]
The operators $K^h_{\bar{x}_0}$ have discrete spectrum (cf, for
instance, \cite{HelNo85,HM88}). Using the simple dilation $X\mapsto
h^{\frac{1}{k+2}}X$, one can show that the operator
$K^h_{\bar{x}_0}$ is unitarily equivalent to
$h^{\frac{2k+2}{k+2}}K^1_{\bar{x}_0}\,$. Thus,
$h^{-\frac{2k+2}{k+2}}K^h_{\bar{x}_0}$ has discrete spectrum,
independent of $h$.

Under the current assumptions, the zero set $U$ of ${\bf B}$ is a finite
collection of points:
\[
U=\{\bar{x}_1,\ldots,\bar{x}_N\}\,.
\]
Let $K^h$ be the self-adjoint operator on $L^2(T_{\bar{x}_1}M)\oplus
\cdots \oplus L^2(T_{\bar{x}_N}M)$ defined by
\begin{equation}\label{e:Kh}
K^h=K^h_{\bar{x}_1}\oplus\cdots\oplus K^h_{\bar{x}_N}\,.
\end{equation}
Let $\mu_0\leq\mu_1\leq \ldots$ be the increasing sequence of
eigenvalues associated with $K^h$ for $h=1$.

\begin{theorem}[\cite{HM}, Theorem 2.5]
For any natural $m$, the eigenvalue $\lambda_m(H^h)$ has an
asymptotic expansion, when $h\to 0$, of the form
\[
\lambda_m(H^h)=h^{(2k+2)/(k+2)}[\mu_m+O(h^{1/(k+2)})]\,.
\]
\end{theorem}

Moreover (\cite[Proposition 2.7]{HM}), if $\mu$ is a non degenerate
eigenvalue of $K^h_{\bar{x}_j}$, for some $j$, then there exists an
eigenvalue $\lambda(H^h)$ of $H^h$ which has a complete asymptotic
expansion of the form
\[
\lambda(H^h)\sim h^{(2k+2)/(k+2)}\sum_{j=0}^{+\infty}a_jh^{j/(k+2)}\,,
\]
with $a_0=\mu$.

\subsection{The case $b_0\neq 0$}
In this subsection, we consider the case when $M$ is a
two-dimensional compact manifold and $b_0\neq 0$. We assume that $M$
has non-empty boundary and the operator $H^h$ satisfies the
Dirichlet boundary conditions. Moreover, we suppose that there is a
unique minimum point $\bar x_0$, which belongs to the interior of
$M$, such that $b(\bar x_0)=b_0$ and which is non degenerate:
\[
{\rm Hess}\,b(\bar x_0)>0.
\]
We introduce in this case the notation
\[
a=\operatorname{Tr} \left(\frac12 {\rm Hess}\,b(\bar
x_0)\right)^{1/2}.
\]

\begin{theorem}[\cite{HM01}\,, Theorem 7.2]
There exist a constant $C>0$ and $h_0>0$, such that, for $h\in (0,h_0]\,$,
\[
-C\,h^{19/8}\leq \lambda_0(H^h)-hb_0-\frac{a^2}{2b_0}h^2\leq C\,h^{5/2}\,.
\]
\end{theorem}

The proof is based on the analysis of the simpler model in $\RR^2$
where near $0$
\[
b(x,y)=b_0+\alpha x^2+\beta y^2\,.
\]
In this case one can also choose a gauge ${\bf
A}=A_1(x,y)dx+A_2(x,y)dy$ such that
\[
A_1(x,y)=0\quad \text{and}\quad A_2(x,y)=b_0+\frac{\alpha}{3}
x^3+\beta xy^2\,.
\]

We mention two open problems in this setting:
\begin{enumerate}
  \item Proof of the existence of a complete asymptotic expansion
  for $\lambda_0(H^h)$ in the two-dimensional case.
  \item Accurate analysis of the bottom of the spectrum in the three-dimensional case.
\end{enumerate}

One should note that the situation is completely different when the
Neumann boundary condition is considered. For a discussion of this case, we refer
the reader to \cite{Fournais-Helffer06} and the references therein.

\section{Hypersurface wells}
In this section, we consider the case when $b_0=0$ and the zero set
$U$ of the magnetic field is a smooth oriented hypersurface $S$.
Moreover, there are constants $k\in \ZZ$, $k>0$, and $C>0$ such
that, for all $x$ in a neighborhood of $S$, we have:
\begin{equation}\label{YK:B1}
C^{-1}\,d(x,S)^k\leq |B(x)|  \leq C\,d(x,S)^k\,.
\end{equation}

This model was introduced for the first time by
Montgomery~\cite{Mont}
 and was further studied in
\cite{HM,Pan-Kwek,syrievienne,Qmath10,miniwells}.

We begin with a discussion of some family of ordinary differential
operators, which play a very important role in the study of this
case.

\subsection{Some ordinary differential operators}\label{s:model} For
any $\alpha\in \RR$ and $\beta\in\RR, \beta\neq 0$, consider the
self-adjoint second order differential operator in $L^2(\RR)$ given
by
\[
Q(\alpha, \beta)=-\frac{d^2 }{d t^2}+ \left(\frac{1}{k+1}\beta
t^{k+1}- \alpha \right)^2.
\]
In the context of magnetic bottles, this family of operators (for
$k=1$) first appears in \cite{Mont} (see also \cite{HM}). Denote by
$\lambda_0(\alpha,\beta)$ the bottom of the spectrum of the operator
$Q(\alpha,\beta)$.

Recall some properties of $\lambda_0(\alpha,\beta)$, which were
established in \cite{Mont,HM,Pan-Kwek}. First of all,
$\lambda_0(\alpha,\beta)$ is a continuous function of $\alpha\in
\RR$ and $\beta\in \RR\setminus\{0\}$. One can see by scaling that,
for $\beta>0$,
\begin{equation}\label{e:scaling}
\lambda_0(\alpha,\beta)= \beta^{\frac{2}{k+2}}\, \lambda_0
(\beta^{-\frac{1}{k+2}}\alpha,1)\,.
\end{equation}
A further discussion depends on $k$ odd or $k$ even.

When $k$ is odd, $\lambda_0(\alpha,1)$ tends to $+\infty$ as $\alpha
\rightarrow -\infty$ by monotonicity. For analyzing its behavior as
$\alpha\rightarrow +\infty$, it is suitable to do a dilation
$t=\alpha^{\frac{1}{k+1}} s$, which leads to the analysis of
\[
\alpha^2 \left(-h^2\frac{d^2}{ds^2} + \left(\frac{s^{k+1}}{k+1}
-1\right)^2 \right)\,,
\]
with $h=\alpha^{-(k+2)/(k+1)}$ small. One can use the semi-classical
analysis (see \cite{CDS} for the one-dimensional case and
\cite{Simon,HSI} for the multidimensional case) to show that
\[
\lambda_0(\alpha,1) \sim (k+1)^{\frac{2k}{k+1}}
\alpha^{\frac{k}{k+1}}\,,\;\mbox{ as } \alpha \rightarrow +\infty\,.
\]
In particular, we see that $\lambda_0(\alpha,1)$ tends to $+\infty\,$.

When $k$ is even, we have $
\lambda_0(\alpha,1)=\lambda_0(-\alpha,1)\,$, and, therefore, it is
sufficient to consider the case $\alpha \geq 0\,$. As $\alpha
\rightarrow +\infty$, semi-classical analysis again shows that
$\lambda_0(\alpha,1)$ tends to $+\infty\,$.

So in both cases, it is clear that the continuous function
$\lambda_0(\alpha,1)$ is positive:
\[
\hat{\nu}:=\inf_{\alpha\in \RR}\lambda_0(\alpha,1)\geq 0\,,
\]
and there exists (at least one) $\alpha_{\mathrm{min}}\in \RR$ such
that $\lambda_0(\alpha,1)$ is minimal:
\[
\lambda_0(\alpha_{\mathrm{min}},1)=\hat{\nu}\,.
\]

The results of numerical computations\footnote{performed for us by
V.~Bonnaillie-No\"el} for $\alpha_{\mathrm{min}}$, $\hat\nu $ and
the second eigenvalue $\lambda_{1}$ of the operator $Q(\alpha_{\rm
min},1)$ are given in Table~\ref{t:comp}.

\begin{table}[h]
\begin{center}
\begin{tabular}{|c|c|c|c|c|c|c|c|}
\hline $k$     &    1  &       2  &     3   &    4   &    5 & 6 & 7\\
\hline $\alpha_{\mathrm{min}}$  & 0.35  &    0    &       0.16 & 0
& 0.10 & 0 &   0.07\\
\hline $\hat\nu $  &       0.57  &  0.66 & 0.68 & 0.76 & 0.81 & 0.87
& 0.92\\ \hline $\lambda_{1}$ & 1.98 & 2.50 & 2.61 & 2.98 & 3.18 &
3.47 & 3.66\\ \hline
\end{tabular}
\end{center}
\caption{Numerical results for $\alpha_{\mathrm{min}}$, $\hat\nu $
and  $\lambda_{1}$}\label{t:comp}
\end{table}

\begin{figure}[h]
\begin{center}
\includegraphics{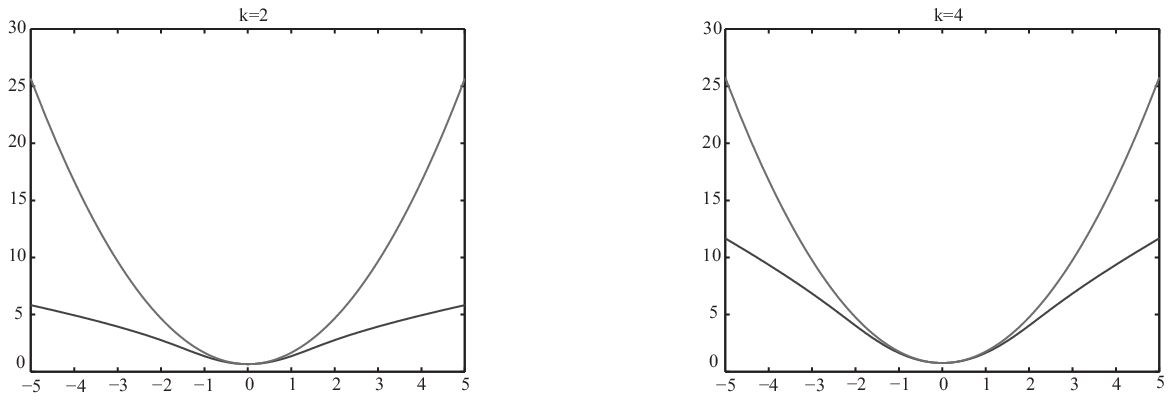}
\end{center}
\caption{$k$ even}\label{f:even}
\end{figure}

\begin{figure}[h]
\begin{center}
\includegraphics{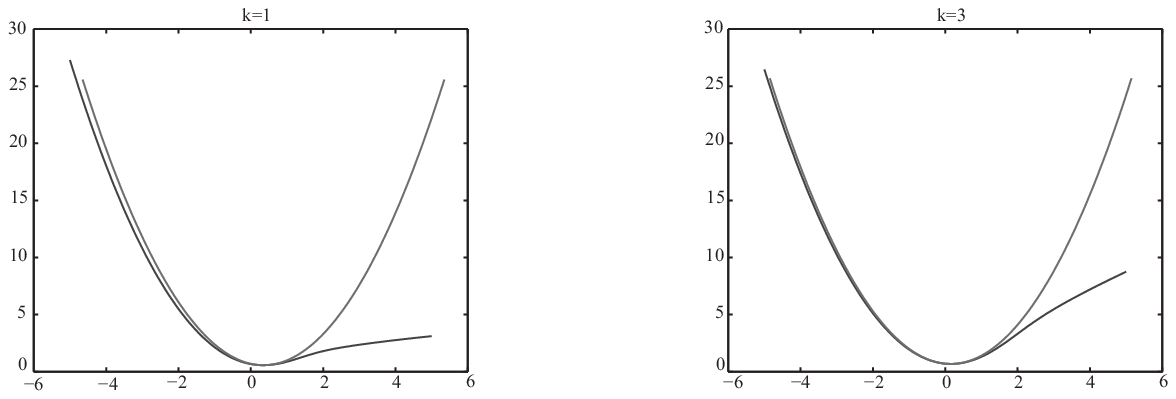}
\end{center}
\caption{$k$ odd} \label{f:odd}
\end{figure}

In Figures~\ref{f:even} and~\ref{f:odd}, one can also see the graphs
of the function $\lambda=\lambda_0(\alpha,1)$ and its quadratic
approximation at $\alpha=\alpha_{\mathrm{min}}$:
\[
\lambda^{\mathrm{quad}}(\alpha)=\lambda_0(\alpha_{\rm
min},1)+\frac12 \frac{\partial^2 \lambda_0}{\partial\alpha^2}
(\alpha_{\mathrm{min}},1)(\alpha- \alpha_{\mathrm{min}} )^2\,.
\]

Numerical computations show that when $k$ is even the minimum is
attained at zero: $\alpha_{\rm min}=0\,$. They also suggest that the
minimum $\alpha_{\rm min}$ is non degenerate:
\[
\frac{\partial^2 \lambda_0}{\partial\alpha^2}
(\alpha_{\mathrm{min}},1) >0\,.
\]
and that the second derivative $\frac{\partial^2
\lambda_0}{\partial\alpha^2} (\alpha_{\mathrm{min}},1)$ tends as $k$
tends to $\infty$ to $2$.

Let $u^0_\alpha\in L^2(\RR)$ be the $L^2$ normalized strictly
positive eigenvector of the operator $Q(\alpha,1)$, corresponding to
the eigenvalue $\lambda_0(\alpha,1)$:
\[
Q(\alpha,1)u^0_\alpha = \lambda_0(\alpha,1) u^0_\alpha, \quad
\|u^0_\alpha\|=1\,.
\]
One can show that $u^0_\alpha$ depends smoothly on $\alpha$. Then
one can show that
\[
\frac{\partial \lambda_0}{\partial\alpha}(\alpha,1)= -2\int
\left(\frac{t^{k+1}}{k+1}- \alpha \right) (u^0_\alpha(t))^2\,dt
\]
and
\[
\frac{\partial^2
\lambda_0}{\partial\alpha^2}(\alpha,1)=2 -4 \int
\frac{t^{k+1}}{k+1}u^0_\alpha(t)\frac{\partial
u^0_\alpha}{\partial\alpha} \,dt\,.
\]
It follows that
\[
\int \left(\frac{t^{k+1}}{k+1}- \alpha_{\mathrm{min}} \right)
(u^0_\alpha(t))^2\,dt=0\,,
\]
and, for $k$ odd, $\alpha_{\mathrm{min}}=\int \frac{t^{k+1}}{k+1}
(u^0_{\alpha_{\mathrm{min}}}(t))^2\,dt>0\,$. It has been claimed that
this minimum is unique for $k=1$ in \cite{Pan-Kwek} and for
arbitrary odd $k$ in \cite{Aramaki05}.

We also have
\[
\left(Q(\alpha,1)-\lambda_0(\alpha,1)\right)\frac{\partial
u^0_\alpha}{\partial\alpha}= \left[2 \left(\frac{t^{k+1}}{k+1}-
\alpha \right) + \frac{\partial \lambda_0}{\partial
\alpha}(\alpha,1)\right] u^0_\alpha\,.
\]

Finally, we mention the following identity (see \cite{Pan-Kwek},
Proposition 3.5 and the formula (3.14)):
\[
\left\|\left(\frac{1}{k+1}t^{k+1}- \alpha_{\rm min}\right)
u^0_{\alpha_{\rm min}}\right\|^2=\frac{\hat\nu}{k+2}\,.
\]

Motivated by numerical computations, we state two conjectures, which
will be very important in further investigations.

\begin{conjecture}\label{c:main}
Any minimum of $\lambda_0(\alpha,1)$ is non-degenerate, that, is,
for any $\alpha_{\rm min}\in \RR$ such that $\lambda_0(\alpha_{\rm
min},1)=\hat{\nu}$ we have
\[
\frac{\partial^2 \lambda_0}{\partial \alpha ^2}(\alpha_{\rm
min},1)>0\,.
\]
\end{conjecture}

\begin{conjecture}\label{c:main2}
There exists a unique $\alpha_{\rm min}\in \RR$ such that
$\lambda_0(\alpha_{\rm min},1)=\hat{\nu}\,$.
\end{conjecture}

One can show that the limit of $\hat \nu$ as $k\rightarrow +\infty$
is $\frac{\pi^2}{4}$, which is the lowest eigenvalue of the Dirichlet
problem for the operator $-d^2/dt^2$ on $(-1,+1)$,
and that Conjecture~\ref{c:main} is true for $k$ large enough.

\subsection{Eigenvalue estimates}
Suppose that the assumption \eqref{YK:B1} holds. Denote by $N$ the
external unit normal vector to $S$ and by $\tilde{N}$ an arbitrary
extension of $N$ to a smooth vector field on $U$. Let $\omega_{0,1}$
be the smooth one form on $S$ defined, for any vector field $V$ on
$S$, by the formula
\[
\langle V,\omega_{0,1}\rangle(y)=\frac{1}{k!}\tilde{N}^k({\mathbf
B}(\tilde{N},\tilde{V}))(y)\,, \quad y\in S\,,
\]
where $\tilde{V}$ is a $C^\infty$ extension of $V$ to $U$. By
\eqref{YK:B1}, it is easy to see that $\omega_{0,1}(x)\not=0$ for
any $x\in S\,$. Denote
\[
\omega_{\mathrm{min}}(B)=\inf_{x\in S} |\omega_{0,1}(x)|>0\,.
\]

As above, $\lambda_0(H^h)$ denotes the bottom of the spectrum of the
operator $H^h$ in $L^2(M)$.

\begin{theorem}[\cite{miniwells}]\label{t:estimate-c}
There exists $C >0$ and $h_0>0$ such that, for any $h\in (0,h_0]$,
we have~:
\[
\hat{\nu}\,
\omega_{\mathrm{min}}(B)^{\frac{2}{k+2}}h^{\frac{2k+2}{k+2}} - C\,
h^{\frac{6k+8}{3(k+2)}} \leq \lambda_0(H^h) \leq \hat{\nu}\,
\omega_{\mathrm{min}}(B)^{\frac{2}{k+2}}h^{\frac{2k+2}{k+2}} + C\,
h^{\frac{6k+8}{3(k+2)}}\,.
\]
\end{theorem}

Observe that a similar result was obtained for the bottom of the
spectrum of the Neumann realization of the operator $H^h$ in a
bounded domain in $\RR^2$ by Pan and Kwek~\cite{Pan-Kwek} in the
case $k=1$ and by Aramaki \cite{Aramaki05} in the case $k$ arbitrary
odd.

As an immediate consequence of Theorems~\ref{t:estimate-c}
and~\ref{YK:hypersurface}, we obtain estimates for the eigenvalues
of the operator $H^h$.

\begin{corollary}[\cite{miniwells}]
For integer $m\geq 0$, we have
\[
\lim_{h\to 0} h^{-\frac{2k+2}{k+2}}\lambda_m(H^h)=\hat{\nu}\,
\omega_{\mathrm{min}}(B)^{\frac{2}{k+2}}\,.
\]
\end{corollary}

The proof of Theorem~\ref{t:estimate-c} is based on reduction to a
second order differential operator $H^{h,0}$ on $\RR\times S$, which
is obtained by expanding the operator $H^h$ near $S$. It is defined
as follows. Let $G$ be the Riemannian metric on $S$ induced by $g$.
Denote by $dx_{G}$ the corresponding Riemannian volume form on $S$.
Let
\[
\omega_{0.0}=i^*_S{\mathbf A}
\]
be the closed one form on $S$ induced by ${\mathbf A}$, where $i_S$
is the embedding of $S$ to $M$. For any $t\in \RR$, let
$P^h_S\left(\omega_{0,0} +\frac{1}{k+1}t^{k+1}\omega_{0,1}\right)$
be a formally self-adjoint operator in $L^2(S, dx_{G})$ defined by
\begin{multline*}
P^h_S\left(\omega_{0,0}+\frac{1}{k+1}t^{k+1}\omega_{0,1}\right)=
\left(ihd+\omega_{0,0}+\frac{1}{k+1}t^{k+1}\omega_{0,1}\right)^*\\
\times\left(ihd+\omega_{0,0}+\frac{1}{k+1}t^{k+1}\omega_{0,1}\right).
\end{multline*}

The operator $H^{h,0}$ is a self-adjoint operator in $L^2(\RR\times
S, dt\,dx_{G})$ defined by the formula
\[
H^{h,0}=-h^2\frac{\partial^2 }{\partial t^2}+
P^h_S\left(\omega_{0,0}+\frac{1}{k+1}t^{k+1}\omega_{0,1}\right).
\]
By Theorem~2.7 of \cite{HM}, the operator $H^{h,0}$ has discrete
spectrum.

Further analysis based on separation of variables leads to spectral
problems for the ordinary differential operator $Q(\alpha, \beta)$
discussed in Subsection~\ref{s:model}. Consider a toy example
considered in \cite{Mont}. Suppose that $n=2$ and the zero set of
$B$ is a connected smooth curve $\gamma$. Let $t\in [0,L)\cong
S^1_L=\RR/L\ZZ$ be the natural parameter along $\gamma$ ($L$ is the
length of $\gamma$). The operator $H^{h,0}$ acts in $L^2(\RR \times
S^1_L)$ by the formula
\begin{equation}\label{ToyMod}
H^{h,0}=-h^2\frac{\partial^2 }{\partial t^2}+ \left(ih\frac{\partial
}{\partial x}+\alpha_1(x)+\frac{1}{(k+1)!}\beta_1(x)
t^{k+1}\right)^2,\quad t\in \RR, \quad x\in S^1_L.
\end{equation}
Choosing an appropriate gauge, without loss of generality, we can
assume that $\alpha_1(x) \equiv \alpha_1 = {\rm const}$. Assume, for
simplicity, that $\beta_1(x) \equiv \beta_1 = {\rm const}$.
Considering Fourier series, we obtain that the operator $H^{h,0}$ is
unitarily equivalent to a direct sum $\bigoplus_{p\in \ZZ} H(a_p)$,
where $$ a_p=a_p(h):=2\pi h p/L-\alpha_1\,,$$ and $H(a)\,, \;a\in \RR\,,$ is an operator
in $L^2(\RR,dt)$ given by
\[
H(a)=-h^2\frac{\partial^2 }{\partial t^2}+
\left(a-\frac{1}{(k+1)!}\beta_1t^{k+1}\right)^2=h^2 Q(h^{-1}a,
h^{-1}\beta_1)\,.
\]
Using \eqref{e:scaling}, we obtain
\[
\inf \sigma(H^{h,0})=\inf_{p\in \ZZ} \sigma(H(a_p))=
\beta_1^{\frac{2}{k+2}} \, h^{\frac{2k+2}{k+2}}\inf_{p\in \ZZ}
\lambda_0 (h^{-\frac{k+1}{k+2}}a_p,1)\,.
\]
We can always find $p_0\in \ZZ$ such that
\[
\left|h^{-\frac{k+1}{k+2}}a_{p_0}-\alpha_{\mathrm{min}}\right|\leq
\frac{2\pi}{L}\, h^{\frac{1}{k+2}}\,.
\]
Therefore, we obtain that
\[
|\inf \sigma(H^{h,0})- \hat\nu\beta_1^{\frac{2}{k+2}}
h^{\frac{2k+2}{k+2}}| \leq C\, h^{\frac{2k+2}{k+2}}
\left|h^{-\frac{k+1}{k+2}}a_{p_0}-\alpha_{\mathrm{min}}\right|^2\leq
C_1\,h^2.
\]
Observe that $\omega_{(0,1)}=\beta_1 dx$ and
$\omega_{\mathrm{min}}=\beta_1$. So we obtain that
\[
\hat\nu\omega_{\mathrm{min}}^{\frac{2}{k+2}}
h^{\frac{2k+2}{k+2}}-C_1\,h^2 \leq \inf \sigma(H^{h,0})\leq
\hat\nu\omega_{\mathrm{min}}^{\frac{2}{k+2}}
h^{\frac{2k+2}{k+2}}+C_1\,h^2\,.
\]
Remark that these estimates are stronger than the estimates of
Theorem~\ref{t:estimate-c}. As observed by Montgomery \cite{Mont}, in this
case, the eigenvalues splitting $\lambda_1-\lambda_0$ between the
second eigenvalue $\lambda_1$ and the lowest eigenvalue $\lambda_0$
of the operator $H^{h,0}$ is $O(h^2)$ and oscillating between this
upper bound and $o(h^2)$.

Moreover, if we admit that $\alpha_{\rm{min}}$ is the unique
critical point of $\lambda_0(\alpha,1)$ (that implies, in
particular, Conjecture~\ref{c:main2}) then, for any $\alpha_1\neq
0$, one can show that there exist $h_0$ and $p_0$ such that,
for any $p$, such that $|p|\geq p_0$ and $\alpha_1\, p>0\,$, there
exists $h_p\in (0,h_0)$ such that $\lim_{p\rightarrow +\infty} h_p
=0$ and the multiplicity
 of the the lowest eigenvalue of $H^{h_p,0}$ is at least $2$.
This is still true if $\alpha_1=0$ and $k$ is odd. On the contrary,
in the case when $k$ is even, if we only admit
Conjecture~\ref{c:main2}, then the multiplicity is $1$.

Let us treat the case when $\alpha_1 >0$. Take an arbitrary $h_0>0$.
Using the asymptotic behavior of $\lambda_0(\alpha,1)$ at $+\infty$
(one can actually prove the monotonicity), we obtain that there
exists $p_0$ such that, for $p\geq p_0$, we have
\[
\lambda_0( h_0^{-\frac{k+1}{k+2}}a_p(h_0),1) < \lambda_0(
h_0^{-\frac{k+1}{k+2}}a_{p+1}(h_0),1)\;.
\]
On the other hand, we observe that, for a given $p$,
\[
\lim_{h\rightarrow 0} h^{-\frac{k+1}{k+2}}a_p
 = -\infty\,.
\]
Using the monotonicity of $\lambda_0(\alpha,1)$ at $-\infty$, we get
\[
\lambda_0( h^{-\frac{k+1}{k+2}}a_p(h),1) > \lambda_0(
h^{-\frac{k+1}{k+2}}a_{p+1}(h),1)\;,
\]
for $h$ small enough. Hence, for $p\geq p_0$, there exists $h_p\in
(0,h_0)$ such that
\[
\lambda_0( h_p^{-\frac{k+1}{k+2}}a_p(h_p),1) = \lambda_0(
h_p^{-\frac{k+1}{k+2}}a_{p+1}(h_p),1)
\]
Since we admit that $\alpha_{\rm{min}}$ is the unique critical point
of $\lambda_0(\alpha,1)$, we immediately get that, for $p\geq p_0$,
\[
\lambda_0( h_p^{-\frac{k+1}{k+2}}a_p(h_p),1)
 = \inf_{q \in \mathbb Z} \lambda_0(
h_p^{-\frac{k+1}{k+2}}a_q(h_p),1)\,,
\]
and
\[
h_p^{-\frac{k+1}{k+2}}a_p(h_p) \leq \alpha_{\mathrm{min}}
 \leq h_p^{-\frac{k+1}{k+2}}a_{p+1}(h_p)\;.
\]
Hence we have, for $p\geq p_0$,
\[
h_p^{-\frac{k+1}{k+2}}a_p(h_p) \leq \alpha_{\mathrm{min}} +
C\left(\frac{2\pi}{L}\right)^2 h_p^{\frac{2}{k+2}}\leq C_1\,,
\]
this shows that $\lim_{p\rightarrow + \infty} h_p =0$.

Like in the case of the Schr\"odinger operator with electric
potential (see \cite{HelSj5}), one can introduce an internal notion
of magnetic well for the fixed hypersurface $S$ in the zero set of
the magnetic field $\mathbf B$. Such magnetic wells can be naturally
called magnetic miniwells. They are defined by means of the function
$|\omega_{0,1}|$ on $S$. Assuming that there exists a non-degenerate
miniwell on $S$, we prove stronger upper bounds for the eigenvalues
of $H^h$.

\begin{theorem}[\cite{miniwells}] \label{t:upperbounds}
Assume that there exist $x_1\in S$ and  $C_1>0$, such that $|\omega_{0,1}(x_1)|=
\omega_{\mathrm{min}}(B)$ and, for all $x\in S$ in some neighborhood
of $x_1$, we have the estimate
\[
C_1^{-1}\,d_S(x,x_1)^2  \leq |\omega_{0,1}(x)|-\omega_{\mathrm{min}}(B)\leq
C_1\,d_S(x,x_1)^2\,.
\]
Then, for any natural $m$, there exist $\widehat C_m>0$ and $h_m>0$ such that,
for any $h\in (0,h_m]\,$, we have
\[
\lambda_m(H^h) \leq \hat{\nu}\,
\omega_{\mathrm{min}}(B)^{\frac{2}{k+2}}\,h^{\frac{2k+2}{k+2}} +
\widehat C_m\,
h^{\frac{2k+3}{k+2}}\,.
\]
\end{theorem}

For the proof of Theorem~\ref{t:upperbounds}, we use a more refined
model operator than the operator $H^{h,0}$ , which is obtained by
considering further terms in the asymptotic expansion of the
operator $H^h$ near $S$. Then we apply the method initiated by
Grushin \cite{Grushin} (and references therein) and Sj\"ostrand
\cite{Sj73} in the context of hypoellipticity.
 We refer also the reader to
\cite{syrievienne} for a discussion of a toy model of this type.

We believe that, if we assume that there exists a unique miniwell
and that Conjecture~\ref{c:main} is true, then, using the methods of
\cite{Fournais-Helffer06}, one can prove the lower bound for the
ground state energy $\lambda_0(H^h)$ of the form
\[
\lambda_0(H^h) \geq \hat{\nu}\,
\omega_{\mathrm{min}}(B)^{\frac{2}{k+2}}h^{\frac{2k+2}{k+2}} - C
h^{\frac{2k+3}{k+2}}\,,
\]
and the upper bound for the splitting between $\lambda_0(H^h)$ and
$\lambda_1(H^h)$ of the form
\[
\lambda_1(H^h)-\lambda_0(H^h)\leq C h^{\frac{2k+3}{k+2}}\,.
\]
Moreover, if, in addition, Conjecture~\ref{c:main2} is true, we
believe that one can prove the lower bound for the splitting between
$\lambda_0(H^h)$ and $\lambda_1(H^h)$ of the form
\[
\lambda_1(H^h)-\lambda_0(H^h)\geq \frac 1C  h^{\frac{2k+3}{k+2}}\,.
\]
 Hence the situation here is quite
different of the case when $n=2$ and $|\omega_{0,1}(x)|$ is constant
along $S$ discussed by Montgomery \cite{Mont} (see the analysis above
of our toy model \eqref{ToyMod}). Remark that the
question about upper and lower bounds for the eigenvalue splitting
$\lambda_1-\lambda_0$ in the Montgomery case is still open.

\section{Periodic operators}\label{s:periodic}
\subsection{The setting of the problem}
In this section, we discuss the case when $M$ is a noncompact
oriented manifold of dimension $n\geq 2$ equipped with a properly
discontinuous action of a finitely generated, discrete group
$\Gamma$ such that $M/\Gamma$ is compact. Suppose that $H^1(M, \RR)
= 0$, i.e. any closed $1$-form on $M$ is exact.

As an example, one can consider the Euclidean space $\RR^n$ equipped
with an action of $\ZZ^n$ by translations or the hyperbolic plane
$\mathbb H$ equipped with an action of the fundamental group of a
compact Riemannian surface of genus $g\geq 2$.

Let $g$ be a $\Gamma$-invariant Riemannian metric and $\bf B$ a
real-valued $\Gamma$-invariant closed 2-form on $M$. Assume that
$\bf B$ is exact and choose a real-valued 1-form $\bf A$ on $M$ such
that $d{\bf A} = \bf B$.

Throughout in this section, we will assume that the magnetic field
has a periodic set of compact potential wells. More precisely, we
assume that there exist a (connected) fundamental domain $\cF$ and a
constant $\epsilon_0>0$ such that
\begin{equation}\label{YK:tr1}
  {\Tr}^+ (B(x)) \geq b_0+\epsilon_0, \quad x\in \partial\cF.
\end{equation}
For any $\epsilon_1 \leq \epsilon_0$, put
\[
U_{\epsilon_1} = \{x\in \cF\,:\, {\Tr}^+ (B(x)) < b_0+ \epsilon_1\}\,.
\]
Thus $U_{\epsilon_1}$ is an open subset of $\cF$ such that
$U_{\epsilon_1}\cap \partial\cF=\emptyset$ and, for $\epsilon_1 <
\epsilon_0$, $\overline{U_{\epsilon_1}}$ is compact and included in
the interior of $\cF$.

We will discuss gaps in the spectrum of the operator $H^h$, which
are located below the top of potential barriers, that is, on the
interval $[0, h(b_0+\epsilon_0)]$. Here by a gap in the spectrum
$\sigma(T)$ of a self-adjoint operator $T$ in a Hilbert space we
understand any connected component of the complement of $\sigma(T)$
in $\RR$, that is, any maximal interval $(a,b)$ such that $(a,b)\cap
\sigma(T) = \emptyset \,$. The problem of existence of gaps in the
spectra of second order periodic differential operators has been
extensively studied recently (some relevant references can be found,
for instance, in \cite{HempelPost02,Qmath10}).

\subsection{Spectral gaps and tunneling effect}
Using the semiclassical analysis of the tunneling effect, it was
shown in \cite{gaps} that the spectrum of the magnetic Schr\"odinger
operator $H^h$ on the interval $[0, h(b_0+\epsilon_0)]$ is localized
in an exponentially small neighborhood of the spectrum of its
Dirichlet realization inside the wells. This result extends to the
periodic setting the result obtained in \cite{HM} in the case of
compact manifolds. It allows us to reduce the investigation of some
gaps in the spectrum of the operator $H^h$ to the study of the
eigenvalue distribution for a ``one-well'' operator and leads us to
suggest a general scheme of a proof of existence of spectral gaps in
\cite{diff2006}. We disregard the analysis of the spectrum in the
above mentioned exponentially small neighborhoods.

For any domain $W$ in $ M$, denote by $H^h_W$ the unbounded
self-adjoint operator in the Hilbert space $L^2(W)$ defined by the
operator $H^h$ in $\overline{W}$ with Dirichlet boundary conditions.
The operator $H^h_W$ is generated by the quadratic form
\[
u\mapsto q^h_W [u] : = \int_W |(ih\,d+{\bf A})u|^2\,dx
\]
with the domain
\[
\Dom (q^h_W) = \{ u\in L^2(W) : (ih\,d+{\bf A})u \in L^2\Omega^1(W),
u\left|_{\partial W}\right.=0 \},
\]
where $L^2\Omega^1(W)$ denotes the Hilbert space of $L^2$
differential $1$-forms on $W$, $dx$ is the Riemannian volume form on
$M$.

Assume now that the operator $H^h$ satisfies the condition of
\eqref{YK:tr1}. Fix $\epsilon_1>0$ and $\epsilon_2>0$ such that
$\epsilon_1 < \epsilon_2 < \epsilon_0$, and consider the operator
$H^h_D$ associated with the domain $D=\overline{U_{\epsilon_2}}$.
The operator $H^h_D$ has discrete spectrum.

\begin{theorem} \label{t:abstract}
Let $N\geq 1$. Suppose that there exist $h_0>0$, $c>0$ and  $M \geq 1$
 and, for each $h\in
(0,h_0]$,   a subset
$\mu_0^h<\mu_1^h<\ldots <\mu_{N}^h$ of an interval $I(h)\subset [0,
h(b_0+\epsilon_1))$ such that~:
\begin{enumerate}
  \item
\begin{equation*}
\begin{split}
& \mu_{j}^h-\mu_{j-1}^h>ch^M, \quad j=1,\ldots,N,\\
& {\rm dist}(\mu_0^h,\partial I(h))>ch^M,\quad {\rm
dist}(\mu_N^h,\partial I(h))>ch^M\,.
\end{split}
\end{equation*}
  \item Each $\mu_j^h, j=0,1,\ldots,N,$ is an approximate eigenvalue of
the operator $H^h_D$: for some $v_j^h\in C^\infty_c(D)$ we have
\begin{equation*}
\|H^h_Dv_j^h-\mu^h_jv_j^h\|=\alpha_j(h)\|v_j^h\|,
\end{equation*}
where $\alpha_j(h)=o(h^M)$ as $h\to 0$.
\end{enumerate}
Then there exists $h_1\in (0, h_0]$ such that, for $h\in (0,h_1]$,
the spectrum of $H^h$ on the interval $I(h)$ has at least $N$
gaps.
\end{theorem}

\subsection{Results on the existence of spectral gaps}
In \cite{diff2006}, we show that, under the assumption
\eqref{YK:tr1}, the spectrum of the operator $H^h$ has gaps (and,
moreover, an arbitrarily large number of gaps) on the interval $[0,
h(b_0+\epsilon_0)]$ in the semiclassical limit $h\to 0\,$. Under
some additional generic assumption, this result was obtained in
\cite{gaps}.

\begin{theorem}
For any natural $N$, there exists $h_0>0$ such that, for any $h\in
(0,h_0]$, the spectrum of $H^h$ in the interval $[0,
h(b_0+\epsilon_0)]$ has at least $N$ gaps.
\end{theorem}

The case when $b_0=0$ and there are regular discrete wells was
considered in \cite{diff2006}.

\begin{theorem}
Suppose that there exist a zero $\bar{x}_0$ of $B$,
$B(\bar{x}_0)=0$,  some integer $k>0$ and a
positive constant $C$ such that, for all $x$ in some neighborhood of
$x_0\,$,  the estimate holds:
\begin{equation}\label{e:discrete}
C^{-1}d(x,x_0)^k\leq {\Tr}^+ (B(x))  \leq C d(x, x_0)^k.
\end{equation}
Then, for any natural $N$, there exist constants $C_N>0$ and $h_N>0$
such that,  for any
$h\in (0,h_N]\,$,  the part of the spectrum of $H^h$ contained in the
interval $[0,C_Nh^{\frac{2k+2}{k+2}}]$  has at least $N$ gaps.
\end{theorem}

A slightly stronger result was shown in \cite{Ko04} under the
assumptions that $b_0=0$ and each zero $\bar{x}_0$ of $B$ satisfies
\eqref{e:discrete}.

\begin{theorem}\label{t:cmp}
Under the current assumptions, there exists an increasing sequence
$\{\mu_m, m\in\NN \}$, satisfying $\mu_m\to\infty$ as $m\to\infty$,
and,  for any $a$ and $b$ satisfying $\mu_m<a<b < \mu_{m+1}$
with some $m$, $h_m>0$ such that, for $h\in (0,h_m]\,$, the interval $[ah^{\frac{2k+2}{k+2}},
bh^{\frac{2k+2}{k+2}}]\,$ does not meet the spectrum of $H^h$. It
  follows that there exists an arbitrarily large
number of gaps in the spectrum of $H^h$ provided the coupling
constant $h$ is sufficiently small.
\end{theorem}

In this case the zero set $U$ in $\cF$ is a finite collection of
points $\{\bar{x}_1,\ldots,\bar{x}_N\}$. Then the sequence $\{\mu_m,
m\in\NN \}$ in Theorem~\ref{t:cmp} is the increasing sequence of
eigenvalues associated with the operator $K^h$ defined in
\eqref{e:Kh}. The proof of Theorem~\ref{t:cmp} is based on abstract
operator-theoretic results obtained in \cite{KMS},

Now suppose that $b_0=0$ and the zero set of the magnetic field is a
smooth oriented hypersurface $S$. Moreover, there are constants
$k\in \ZZ, k>0$ and $C>0$ such that for all $x\in U$ we have:
\[
C^{-1}d(x,S)^k\leq |B(x)|\leq C d(x,S)^k\,.
\]

First of all, note, that the estimates of Theorem~\ref{t:estimate-c}
hold in this setting \cite{miniwells}. In \cite{Qmath10} we have
proved the following result.

\begin{theorem}\label{YK:hypersurface}
For any $a$ and $b$ such that
\[
\hat{\nu}\, \omega_{\mathrm{min}}(B)^{\frac{2}{k+2}}<a < b\,,
\]
and for any natural $N$, there exists $h_0>0$ such that, for any
$h\in (0,h_0]\,$, the spectrum of $H^h$ in the interval
$
[h^{\frac{2k+2}{k+2}}a, h^{\frac{2k+2}{k+2}}b]
$
has at least $N$ gaps.
\end{theorem}

Finally, assuming the existence of a non-degenerate miniwell on $S$,
we prove the existence of gaps in the spectrum of $H^h$ on intervals
of size $h^{\frac{2k+3}{k+2}}$, close to the bottom
$\lambda_0(H^h)$.

\begin{theorem}\label{t:gaps}
Under the current assumptions, suppose that there exist $x_1\in S$
and  $C_1>0$, such
that $|\omega_{0,1}(x_1)|= \omega_{\mathrm{min}}(B)$ and, for all
$x\in S$ in some neighborhood of $x_1$
\[
\frac{1}{C_1}\,d_S(x,x_1)^2  \leq  |\omega_{0,1}(x)|-\omega_{\mathrm{min}}(B)\leq
C_1\, d_S(x,x_1)^2.
\]
Then, for any natural $N$, there exist $b_N>0$ and $h_N>0$ such that,
  for any $h\in (0,h_N]\,$,
the spectrum of $H^h$ in the interval
\[
\left[\hat{\nu}\,
\omega_{\mathrm{min}}(B)^{\frac{2}{k+2}}\,h^{\frac{2k+2}{k+2}}\,,\,
\hat{\nu}\,
\omega_{\mathrm{min}}(B)^{\frac{2}{k+2}}\,h^{\frac{2k+2}{k+2}} +
b_N\,h^{\frac{2k+3}{k+2}}\right]
\]
has at least $N$ gaps.
\end{theorem}


\begin{thebibliography}{10}

\bibitem{Aramaki05}
Junichi Aramaki.
\newblock Asymptotics of upper critical field of a superconductor in applied
  magnetic field vanishing of higher order.
\newblock {\em Int. J. Pure Appl. Math.}, 21(2):151--166, 2005.

\bibitem{CDS}
J.-M. Combes, P.~Duclos, and R.~Seiler.
\newblock Kre\u\i n's formula and one-dimensional multiple-well.
\newblock {\em J. Funct. Anal.}, 52(2):257--301, 1983.

\bibitem{CFKS}
H.~L. Cycon, R.~G. Froese, W.~Kirsch, and B.~Simon.
\newblock {\em Schr\"odinger operators with application to quantum mechanics
  and global geometry}.
\newblock Texts and Monographs in Physics. Springer-Verlag, Berlin, study
  edition, 1987.

\bibitem{Fournais-Helffer06}
Soeren Fournais and Bernard Helffer.
\newblock Accurate eigenvalue asymptotics for the magnetic {N}eumann
  {L}aplacian.
\newblock {\em Ann. Inst. Fourier (Grenoble)}, 56(1):1--67, 2006.

\bibitem{Fournais-Helffer:book}
Soeren Fournais and Bernard Helffer.
\newblock {\em Spectral methods in surface superconductivity}.
\newblock Birkh\"auser, Basel, 2009.

\bibitem{Grushin}
V.~V. Gru{\v{s}}in.
\newblock Hypoelliptic differential equations and pseudodifferential operators
  with operator-valued symbols.
\newblock {\em Mat. Sb. (N.S.)}, 88(130):504--521, 1972.

\bibitem{Helffer-Kuroda}
Bernard Helffer.
\newblock On spectral theory for {S}chr\"odinger operators with magnetic
  potentials.
\newblock In {\em Spectral and scattering theory and applications}, volume~23
  of {\em Adv. Stud. Pure Math.}, pages 113--141. Math. Soc. Japan, Tokyo,
  1994.

\bibitem{Helffer-ECM}
Bernard Helffer.
\newblock Analysis of the bottom of the spectrum of {S}chr\"odinger operators
  with magnetic potentials and applications.
\newblock In {\em European {C}ongress of {M}athematics}, pages 597--617. Eur.
  Math. Soc., Z\"urich, 2005.

\bibitem{syrievienne}
Bernard Helffer.
\newblock Introduction to semi-classical methods for the schr\"odinger operator
  with magnetic fields.
\newblock In {\em Aspects th\'eoriques et appliqu\'es de quelques EDP issues de
  la g\'eom\'etrie ou de la physique, Proceedings of the CIMPA School held in
  Damas (Syrie) (2004)}, S\'eminaires et Congr\`es. SMF, 2009.

\bibitem{diff2006}
Bernard Helffer and Yuri~A. Kordyukov.
\newblock The periodic magnetic {Schr\"odinger} operators: spectral gaps and
  tunneling effect.
\newblock {\em Trudy Matematicheskogo Instituta Imeni V.A. Steklova},
  261:176--187, 2008.
\newblock translation in Proceedings of the Steklov Institute of Mathematics,
  261 (2008) 171-182.

\bibitem{gaps}
Bernard Helffer and Yuri~A. Kordyukov.
\newblock Semiclassical asymptotics and gaps in the spectra of periodic
  {S}chr\"odinger operators with magnetic wells.
\newblock {\em Trans. Amer. Math. Soc.}, 360(3):1681--1694 (electronic), 2008.

\bibitem{Qmath10}
Bernard Helffer and Yuri~A. Kordyukov.
\newblock Spectral gaps for periodic {Schr\"odinger} operators with
  hypersurface magnetic wells.
\newblock In {\em ``Mathematical results in quantum mechanics'', Proceedings of
  the QMath10 Conference Moieciu, Romania 10 - 15 September 2007}, pages
  137--154. World Sci. Publ., Singapore, 2008.

\bibitem{miniwells}
Bernard Helffer and Yuri~A. Kordyukov.
\newblock Spectral gaps for periodic {Schr\"odinger} operators with
  hypersurface magnetic wells: {Analysis} near the bottom.
\newblock E-print arXiv:0812.4350 [math.SP], 2008.

\bibitem{HM}
Bernard Helffer and Abderemane Mohamed.
\newblock Semiclassical analysis for the ground state energy of a
  {S}chr\"odinger operator with magnetic wells.
\newblock {\em J. Funct. Anal.}, 138(1):40--81, 1996.

\bibitem{HM88}
Bernard Helffer and Abderemane Mohamed.
\newblock Asymptotic of the density of states for the {S}chr\"odinger operator
  with periodic electric potential.
\newblock {\em Duke Math. J.}, 92(1):1--60, 1998.

\bibitem{HM01}
Bernard Helffer and Abderemane Morame.
\newblock Magnetic bottles in connection with superconductivity.
\newblock {\em J. Funct. Anal.}, 185(2):604--680, 2001.

\bibitem{HelNo85}
Bernard Helffer and Jean Nourrigat.
\newblock {\em Hypoellipticit\'e maximale pour des op\'erateurs polyn\^omes de
  champs de vecteurs}, volume~58 of {\em Progress in Mathematics}.
\newblock Birkh\"auser Boston Inc., Boston, MA, 1985.

\bibitem{Helffer-Robert84}
Bernard Helffer and Didier Robert.
\newblock Puits de potentiel g\'en\'eralis\'es et asymptotique semi-classique.
\newblock {\em Ann. Inst. H. Poincar\'e Phys. Th\'eor.}, 41(3):291--331, 1984.

\bibitem{HSI}
Bernard Helffer and Johannes Sj{\"o}strand.
\newblock Multiple wells in the semiclassical limit. {I}.
\newblock {\em Comm. Partial Differential Equations}, 9(4):337--408, 1984.

\bibitem{HelSj5}
Bernard Helffer and Johannes Sj{\"o}strand.
\newblock Puits multiples en m\'ecanique semi-classique. {V}. \'{E}tude des
  minipuits.
\newblock In {\em Current topics in partial differential equations}, pages
  133--186. Kinokuniya, Tokyo, 1986.

\bibitem{HempelPost02}
Rainer Hempel and Olaf Post.
\newblock Spectral gaps for periodic elliptic operators with high contrast: an
  overview.
\newblock In {\em Progress in analysis, {V}ol.\ {I}, {II} ({B}erlin, 2001)},
  pages 577--587. World Sci. Publ., River Edge, NJ, 2003.

\bibitem{Ko04}
Yuri~A. Kordyukov.
\newblock Spectral gaps for periodic {S}chr\"odinger operators with strong
  magnetic fields.
\newblock {\em Comm. Math. Phys.}, 253(2):371--384, 2005.

\bibitem{KMS}
Yuri~A. Kordyukov, Varghese Mathai, and Mikhail Shubin.
\newblock Equivalence of spectral projections in semiclassical limit and a
  vanishing theorem for higher traces in {$K$}-theory.
\newblock {\em J. Reine Angew. Math.}, 581:193--236, 2005.

\bibitem{Mall}
Paul Malliavin.
\newblock Analyticit\'e transverse d'op\'erateurs hypoelliptiques {$C\sp 3$}
  sur des fibr\'es principaux. {S}pectre \'equivariant et courbure.
\newblock {\em C. R. Acad. Sci. Paris S\'er. I Math.}, 301(16):767--770, 1985.

\bibitem{Mohamed-Raikov}
Abd{\'e}r{\'e}mane Mohamed and George~D. Ra{\u\i}kov.
\newblock On the spectral theory of the {S}chr\"odinger operator with
  electromagnetic potential.
\newblock In {\em Pseudo-differential calculus and mathematical physics},
  volume~5 of {\em Math. Top.}, pages 298--390. Akademie Verlag, Berlin, 1994.

\bibitem{Mont}
Richard Montgomery.
\newblock Hearing the zero locus of a magnetic field.
\newblock {\em Comm. Math. Phys.}, 168(3):651--675, 1995.

\bibitem{Pan-Kwek}
Xing-Bin Pan and Keng-Huat Kwek.
\newblock Schr\"odinger operators with non-degenerately vanishing magnetic
  fields in bounded domains.
\newblock {\em Trans. Amer. Math. Soc.}, 354(10):4201--4227 (electronic), 2002.

\bibitem{Simon}
Barry Simon.
\newblock Semiclassical analysis of low lying eigenvalues. {I}. {N}ondegenerate
  minima: asymptotic expansions.
\newblock {\em Ann. Inst. H. Poincar\'e Sect. A (N.S.)}, 38(3):295--308, 1983.

\bibitem{Sj73}
Johannes Sj{\"o}strand.
\newblock Operators of principal type with interior boundary conditions.
\newblock {\em Acta Math.}, 130:1--51, 1973.

\bibitem{Ue}
Naomasa Ueki.
\newblock Lower bounds for the spectra of {S}chr\"odinger operators with
  magnetic fields.
\newblock {\em J. Funct. Anal.}, 120(2):344--379, 1994.

\end{thebibliography}
\end{document}